\documentclass[12pt]{amsart}

\usepackage{color}
\usepackage{amsbsy} 
\usepackage{amssymb,amsmath,amstext,mathrsfs}
 
\newtheorem{theorem}{Theorem}[section]
\newtheorem{lemma}[theorem]{Lemma}
\newtheorem{proposition}[theorem]{Proposition}
\newtheorem{corollary}[theorem]{Corollary}
\theoremstyle{definition}
\newtheorem{definition}[theorem]{Definition}
\newtheorem{example}[theorem]{Example}
\newtheorem{remark}[theorem]{Remark}

\newtheorem{step}{Step}

\numberwithin{equation}{section}

\newcommand{\va}{\left|}
\newcommand{\vb}{\right|}
\newcommand{\vd}{\right\| }
\newcommand{\vc}{\left\| }

\newcommand{\pl}{\left(}
\newcommand{\pr}{\right)}

\newcommand{\tl}{\left\{}
\newcommand{\tr}{\right\}}

\newcommand{\ra}{\rightarrow}

\newcommand{\cx}{\mathrm{Lip} (X)}

\newcommand{\cxe}{\mathrm{Lip} (X, E)}
\newcommand{\cyf}{\mathrm{Lip} (Y, F)}
\newcommand{\cxea}{\mathrm{Lip}^{\alpha} (X, E)}
\newcommand{\cyfa}{\mathrm{Lip}^{\alpha} (Y, F)}
\newcommand{\cxeh}{\mathrm{Lip} \pl \widehat{X}, E \pr}
\newcommand{\cyfh}{\mathrm{Lip} \pl \widehat{Y}, F \pr}

\newcommand{\axe}{\mathrm{Lip}_{\ah} (X, E)}
\newcommand{\bxe}{\mathrm{Lip}_{\bh} (X, E)}

\newcommand{\aie}{\mathrm{Lip} (\mathfrak{A} \pl T^{-1} \pr , E)}
\newcommand{\bie}{\mathrm{Lip} (\mathfrak{B} \pl T^{-1} \pr , E)}
\newcommand{\aif}{\mathrm{Lip} (\mathfrak{A} (T)  , F)}
\newcommand{\bif}{\mathrm{Lip} (\mathfrak{B} (T), F)}

\newcommand{\as}{\mathfrak{A} \pl T^{-1} \pr }
\newcommand{\bs}{\mathfrak{B} \pl T^{-1} \pr }

\newcommand{\at}{\mathfrak{A} (T)}
\newcommand{\bt}{\mathfrak{B} (T)}

\newcommand{\ash}{\mathfrak{A} \pl T_{\ah}^{-1} \pr }

\newcommand{\ath}{\mathfrak{A} (T_{\ah})}
\newcommand{\bth}{\mathfrak{B} (T_{\ah})}

\newcommand{\ah}{\mathfrak{A}}
\newcommand{\bh}{\mathfrak{B}}

\newcommand{\iso}{\mathrm{iso}_{<2} (Y,X)}

\newcommand{\isol}{\mathrm{iso}_{<2} \pl \widehat{Y}, \widehat{X} \pr }

\newcommand{\isob}{\mathrm{iso}_{<2} \pl \bt, \bs \pr }

\newcommand{\isoa}{\mathrm{iso}_{<2} \pl \at, \as \pr }

 \newcommand{\ha}{h_{\ah}}

 \newcommand{\hb}{h_{\bh}}
 
 \newcommand{\ja}{J_{\ah}}
 
\newcommand{\jb}{J_{\bh}}

\begin{document}

\title[noncompactness and noncompleteness]{Noncompactness and noncompleteness in isometries of Lipschitz  spaces}

\author{Jes\'us Araujo and Luis Dubarbie}

\address{Departamento de Matem\'aticas, Estad\'istica y
Computaci\'on, Facultad de Ciencias, Universidad de Cantabria,
Avenida de los Castros s/n, E-39071, Santander, Spain.}

\email{araujoj@unican.es} \email{luis.dubarbie@gmail.com}

\thanks{Research partially supported by the Spanish Ministry of Science and Education
(MTM2006-14786). The second author was partially supported by a predoctoral
grant from the University of Cantabria and the Government of Cantabria.}



\begin{abstract}
We solve the following  three  questions concerning surjective linear isometries between spaces
of Lipschitz functions $\mathrm{Lip}(X,E)$ and $\mathrm{Lip}(Y,F)$, for strictly convex normed spaces $E$ and $F$ and
metric spaces $X$ and $Y$:
\begin{enumerate} 
\item Characterize those base spaces $X$ and $Y$ for which all isometries are weighted composition maps.
\item Give a condition independent of base spaces under which all isometries are weighted composition maps.
\item Provide the general form of an isometry, both when it is a weighted composition map and when it is not.
\end{enumerate}
In particular, we prove that requirements of completeness on $X$ and $Y$ are   not
necessary when $E$ and $F$ are not complete, which is in sharp contrast with results known in the scalar context. 
\end{abstract}

\maketitle

\section{Introduction}

It is well known that not all surjective (linear) isometries between spaces of Lipschitz functions
on general metric spaces $X$ and $Y$ can be written as weighted composition maps (see for instance  \cite[p. 61]{W}). 
Attempts to identify the isometries which can be described in that way have been  done in three ways,
each  trying to provide an answer to one of the following  questions:
\begin{enumerate} 
\item\label{rabanito} Characterize those base spaces $X$ and $Y$ for which all isometries are weighted composition maps.
\item\label{carolo} Give a condition independent of base spaces under which all isometries are weighted composition maps.
\item\label{ocumier} Provide the general form of an isometry, both when it is a weighted composition map and when it is not.
\end{enumerate}

The first question was 
studied by Weaver  for a general metric  in the {\em scalar-valued} setting 
(see \cite{W1} or \cite[Section 2.6]{W}), 
and the second one has been 
recently  treated  by Jim\'enez-Vargas and  Villegas-Vallecillos in the more general setting of
{\em vector-valued} functions (see \cite{JV1}). 
In the latter,  the Banach spaces   where the functions take values are assumed to be strictly convex. 
This is certainly not a heavy restriction, as this type of results is known not to hold
for general Banach spaces. Strict convexity is actually a very common and reasonable assumption, 
even if, at least in other contexts, it is not the unique possible (see for instance \cite{A5, B, FJ1, Je}).
 As for the third question, an answer was given by Mayer-Wolf for compact base spaces in the scalar context, not for a general metric $d$, but for powers  $d^{\alpha}$
with $0< \alpha <1$.

 Weaver
 proved that completeness and $1$-connectedness of 
 $X$ and $Y$ are
sufficient conditions, and that the weighted composition isometries must have  a very special form. 
More concretely, given complete $1$-connected metric spaces $X$ and $Y$ with diameter at most $2$, a linear bijection 
$T :\mathrm{Lip}(X) \ra \mathrm{Lip}(Y)$ is an isometry if and only if $T f = \alpha \cdot f \circ h$
 for every $f$, where
$\alpha \in \mathbb{K}$, $\va \alpha \vb =1$, and $h:Y \ra X$ is an isometry. Requirements of $1$-connectedness on both $X$ and $Y$  (that is, they cannot be decomposed into two nonempty
disjoint sets whose distance is greater than or equal to $1$) cannot be dropped  in general.
And, obviously, $\mathrm{Lip}(X)$ and $\mathrm{Lip}(\widehat{X})$ are linearly
isometric when   $X$ is not complete (where $\widehat{X}$ denotes the completion of $X$), so requirements of completeness cannot
be dropped either.

On the other hand, Jim\'enez-Vargas and  Villegas-Vallecillos 
gave
a general representation  in the spirit of
the classical Banach-Stone Theorem (along with related results for isometries not necessarily surjective). 
Assumptions  include compactness of base metric spaces
and the fact that the isometry fixes
a (nonzero) constant function. The conclusion in the surjective case is that the isometry $T: \mathrm{Lip}(X,E) \ra \mathrm{Lip}(Y,E)$ is of the form  $ T f( y) = J y (f ( h (y)))$  for all 
$ f\in {\rm Lip}(X,E)  $ and $y\in Y$,
where $h$ is a bi-Lipschitz homeomorphism  (that is, $h$ and $h^{-1}$ are Lipschitz) from $Y$ onto $X$ 
and $J$ is a 
Lipschitz map   from $Y$ into the set  $I(E,E)$
 of all surjective linear isometries on  the strictly convex Banach space  $E$.
They also proved that this result can be sharpened under stronger hypotheses,
 but the above 
assumptions remain basically the same, so that the results do not
provide an "if and only if" description.

Finally, Mayer-Wolf not only characterized the family of compact spaces for which the associated
Lipschitz spaces admitted isometries that were not composition operators, but also gave their general form. In principle, it is not clear  whether or not his results   can be extended  to spaces endowed with a metric not 
 of the form $d^{\alpha}$. In fact, the answer, as we will see here, is not completely positive.

The aim of this paper is to give, in the {\em vector-valued} setting, a complete answer to  
questions (\ref{rabanito}),(\ref{carolo})  and (\ref{ocumier})  (just assuming strict convexity of $E$ and $F$). 
The general answer is not known even in the scalar setting, which can be included here as a special case. 
We also prove, on the one hand,  that conditions of compactness can be replaced with just 
completeness on base spaces
and, on the other hand, that
even  completeness can be dropped when the
normed spaces $E$ and $F$ are  not complete
(which is in sharp contrast
with the behaviour in the scalar case). 

To solve (\ref{carolo}), we show that the condition on the preservation of a constant function (as given in \cite{JV1})
can be replaced with a milder one (see Theorem~\ref{marthahoyo}). 
We use it to solve (\ref{ocumier}) (see Theorem~\ref{nonelk}, and more in general 
Theorem~\ref{marthahoyo}, Corollary~\ref{bocadillo} and Remark~\ref{sinpase}). Our answer also applies to
results on metrics $d^{\alpha}$ in \cite{M}, and a key to understand the generalization is 
Proposition~\ref{pesafrank}.  
An answer to (\ref{rabanito}) is given as a direct consequence of the results concerning 
(\ref{carolo})  and (\ref{ocumier})
(see Corollaries~\ref{deluxion} and~\ref{rentreeflot}).  
As a special case we provide   the natural
 counterpart of the description given in \cite{W1} 
 (see Corollary~\ref{siegediesrefe} and Remark~\ref{guelgafunc}). 
 We finally mention that we do not use the same techniques
as in  \cite{W1} nor as  in \cite{JV1}; 
instead we study {\em surjective linear} isometries through biseparating maps, which has proven successful in various contexts (see for instance
\cite{A5, HBN} for recent references).

Other papers where related operators have been recently studied in similar contexts are \cite{AEV, GJ, GJ2, JV, Le} (see also  \cite{FJ, JP, L, M, R, V}).

\section{Preliminaries and notation}

Recall that, given metric spaces $(X, d_1)$ and $(Y, d_2)$, a map $f:X \rightarrow Y$ is said to 
be \emph{Lipschitz} if there
exists a constant $k\geq 0$ such that
$
d_2 ( f(x) , f(y) ) \leq k \ d_1(x,y)
$
for each $x,y\in X$, and that 
 the \emph{Lipschitz number} of $f$ is   
$$L(f):=\mathrm{sup}\left\{\frac{ d_2 (f(x) , f(y) )}{d_1 (x,y)}:x,y\in
X, \ x\neq y\right\}.$$

Given a normed space $E$  (over $\mathbb{K} = \mathbb{R}$ or $\mathbb{C}$), we
denote by  $\mathrm{Lip}(X,E)$ the space of all \emph{bounded} $E$-valued
Lipschitz functions  on $X$. 
We endow   $\mathrm{Lip}(X,E)$ 
 with the norm 
$
\vc \cdot \vd_{L} :=\mathrm{max}\left\{\vc \cdot \vd_{\infty},L(\cdot)\right\}
$
(where $\vc  \cdot \vd_{\infty}$ denotes
the usual supremum norm).
 
As a particular case, we can consider in $X$  a power ${d_1}^{\alpha}$     of   the metric $d_1$,
  $0< \alpha <1$.  The corresponding space of all  bounded  $E$-valued
Lipschitz functions  on $X$ with respect to ${d_1}^{\alpha}$ is then denoted by $\mathrm{Lip}^{\alpha} (X,E)$.

Recall also that  a  normed space $E$ is said to be \emph{strictly convex} if 
$\vc \mathrm{e}_{1} + \mathrm{e}_{2} \vd <2$ whenever $\mathrm{e}_{1},\mathrm{e}_{2}$ are different vectors of 
norm $1$ in $E$ or, equivalently, that 
$\vc \mathrm{e}_{1} + \mathrm{e}_{2} \vd = \vc \mathrm{e}_{1} \vd + \vc \mathrm{e}_{2} \vd $ 
($\mathrm{e}_{1}, \mathrm{e}_{2} \neq 0$)  implies $\mathrm{e}_{1} = \alpha \mathrm{e}_{2}$ 
for some $\alpha >0$
(see \cite[pp. 332--336]{K}). 
From this, it follows that, 
given $\mathrm{e}_{1},\mathrm{e}_{2}\in E\setminus \{0\}$,
\begin{equation}\label{strictly}
\vc \mathrm{e}_{1} \vd , \vc \mathrm{e}_{2} \vd<\max\{\vc \mathrm{e}_{1}+\mathrm{e}_{2} \vd,\vc \mathrm{e}_{1}-\mathrm{e}_{2} \vd\}  ,
\end{equation}
which is an inequality we will often use. The fact that a normed space is strictly convex does not imply
that its completion is.  Indeed every infinite-dimensional separable Banach space can be renormed to be not
stricly convex and to contain a strictly convex dense subspace of codimension one (see \cite{GZ}). 

From now on, unless otherwise stated, we assume that 
$E$ and $F$ are strictly convex normed spaces (including the cases $E=\mathbb{K}$, $F = \mathbb{K}$).

As we mentioned above, on our way to Theorem~\ref{marthahoyo} we will deal with biseparating maps.    Recall that   separating maps are those preserving disjointness of {\em cozero sets} (where the cozero set of
a function 
 $f: X \ra E$ is defined as $c(f):=\{x\in X:f(x)\neq 0\}$). More concretely,  
we will say that
a   linear map $T:\mathrm{Lip}(X,E)\rightarrow \mathrm{Lip}(Y,F)$ is \emph{separating} if 
  $c(Tf)\cap c(Tg)=\emptyset$ whenever $f,g\in \mathrm{Lip}(X,E)$ satisfy $c(f)\cap c(g)=\emptyset$.
Moreover, $T$ is said to be \emph{biseparating} if it is bijective and both $T$ and its inverse
are separating maps.

Obviously, if $f:X \ra E$ is Lipschitz and bounded,
then so  is the map $\vc f \vd : X \ra \mathbb{R}$ defined by $\vc f \vd (x) := \vc f(x) \vd$ for every $x \in X$.
 It is also clear that $\vc f \vd$ can be continuously extended to a Lipschitz function
$\widehat{\vc f \vd} : \widehat{X} \ra \mathbb{R}$ 
defined on  the completion $\widehat{X}$ of $X$. 
More in general, if $x \in \widehat{X} \setminus X$,
we say that $f$ {\em admits an extension to $x$} if it can be continuously extended to a map
$\widehat{f} : X \cup \tl x \tr \ra E$. Clearly, when $E$ is complete and $X$ is not, $f$ admits
a continuous extension to the whole $\widehat{X}$, and the extension $\widehat{f}: \widehat{X} \ra E$
is a Lipschitz function with $\vc \widehat{f} \vd_{\infty} = \vc f \vd_{\infty}$ 
and $L \pl \widehat{f} \pr = L \pl f \pr$.
For this reason, when $E$ and $F$ 
are complete, every surjective linear isometry 
$T:\mathrm{Lip}(X,E)\rightarrow \mathrm{Lip}(Y,F)$ can be associated in a canonical way to
another one $\widehat{T} :\mathrm{Lip} \pl \widehat{X},E \pr \rightarrow \mathrm{Lip} \pl \widehat{Y},F \pr $ 
(which coincides with $T$ only if $X$ and $Y$ are complete).

Given $R>0$, we  define in $X$ the following equivalence relation: we put $x \sim_R y$ if there exist
$x_1, \ldots , x_n \in X$ with $x=x_1$, $y=x_n$, and $d(x_i, x_{i+1} ) <R$ for $i=1, \ldots, n-1$.
We call $R$-component each of the equivalence classes of $X$ by $\sim_R$. The set of all $R$-components in $X$
is denoted by $\mathtt{Comp}_R (X)$.

We say that a bijective map $h: Y \ra X$ preserves distances less than $2$ 
 if $d_1 (h(y),h(y')) = d_2(y,y')$ 
 whenever $d_2(y,y') <2 $. We denote by $\iso$ the set of all maps $h: Y \ra X$ such that both $h$ and $h^{-1} $ 
preserve distances less than $2$. Notice that every $h \in \iso$ is  a homeomorphism  and that, when $X$ is bounded,  then  it is also a Lipschitz map (see also Remark~\ref{subidahoyo}).

\begin{definition}\label{martha}
Let $I(E,F)$ be the set of  all linear isometries from $E$ onto $F$.
 We say that a map $T:\mathrm{Lip}(X,E)\rightarrow \mathrm{Lip}(Y,F)$
is a {\em standard} isometry if there exist  $h \in \iso$  
  and a map $J: Y \ra I(E,F)$ constant on each $2$-component of $Y$ such that 
$$
Tf(y)= Jy (f(h(y)))
$$
for all $f\in \mathrm{Lip}(X,E)$ and $y\in Y$.
\end{definition}

\begin{remark}\label{scientiaedies}
Notice that a standard isometry is indeed a surjective linear isometry. Theorem~\ref{marthahoyo} 
 gives a condition under which both classes of operators coincide. 
Also, when $Y$ is $2$-connected, the map $J$ is constant, so there exists 
a surjective linear isometry  $\mathbf{J} : E \ra F$   such that
$$
Tf(y)= \mathbf{J} (f(h(y)))
$$
for all $f\in \mathrm{Lip}(X,E)$ and $y\in Y$.
In particular, Corollary~\ref{siegediesrefe} roughly says  that this is the only way to obtain  an isometry when  one of the base spaces is 
$1$-connected (see also Remark~\ref{guelgafunc}).
\end{remark}

In the definition of standard isometry, we see that $X$ and $Y$ are very much related. In particular, one is complete if and only if the other is. There  are interesting cases which are {\em almost} standard in some sense. For instance, when $E$ is complete, the natural inclusion
$\mathbf{i}_X : \mathrm{Lip}(X,E) \rightarrow \mathrm{Lip} \pl \widehat{X},E \pr $ is not standard if $X$ is not complete, but we immediately obtain a standard isometry from it in a natural way.

On the other hand, when $E$ and $F$ are   complete,  every $T$ can be 
written as $T = {\mathbf{i}_Y}^{-1} \circ  \widehat{T} \circ \mathbf{i}_X$. In the following definition, we 
distinguish between this kind of isometries 
and 
{\em nonstandard} isometries.

\begin{definition}\label{desk}
We say that a surjective linear isometry $T:\mathrm{Lip}(X,E)\rightarrow \mathrm{Lip}(Y,F)$ is 
{\em nonstandard} if $T$ and $\widehat{T}$ (if it can be defined, that is, if $E$ and $F$ are complete) 
are not standard.
\end{definition}

 Just a special family of spaces allows   defining properly nonstandard isometries. We call them spaces of
 {\em type} $\mathbf{A}$.

\begin{definition}\label{nikolabit}
We say that a metric space $X$ is of   type $\mathbf{A}$ if there are a partition of $X$ into two subsets 
$\mathfrak{A}$, $\mathfrak{B}$, and a   map $\varphi : \mathfrak{A} \ra \mathfrak{B}$ with the following properties: 
\begin{enumerate}
\item\label{cand} $d (x, z) = 1 + d(\varphi(x), z)$ whenever $x \in \mathfrak{A}$ and  $z \in \mathfrak{B}$ satisfy $d(\varphi(x),z) < 1$,
\item\label{ela} $d (x, z) \ge 2$ whenever $x \in \mathfrak{A}$ and  $z \in \mathfrak{B}$  satisfy $d(\varphi(x),z) \ge 1$,
\item $d (x_1, x_2 ) \ge 2 $ whenever $x_1, x_2 \in \mathfrak{A}$ and $\varphi(x_1) \neq \varphi (x_2)$.
\end{enumerate}
For each $E$, the operator  $S_{\varphi} : \cxe \ra \cxe$
defined, for each $f \in \cxe$,  by
\begin{displaymath}
S_{\varphi} f (x) := \left\{ \begin{array}{ll} f (x) & \mbox{if }  x \in \mathfrak{B}
\\
f(\varphi (x)) - f(x) & \mbox{if } x \in \mathfrak{A}
\end{array} \right.
\end{displaymath}
is said to be the {\em purely nonstandard} map associated to $\varphi$.
\end{definition}

\begin{remark}\label{gemacurruscos}
It is easy to check that $S_{\varphi}$ is   linear and bijective. Also
$\vc S_{\varphi} (f) \vd \le 1$ whenever $\vc f \vd \le 1$. Taking into account that
${S_{\varphi}}^{-1} = S_{\varphi}$, this implies that $S_{\varphi}$ is indeed a nonstandard isometry.
Theorem~\ref{nonelk} and Remark~\ref{sinpase} basically say that every nonstandard isometry
is the composition of a standard and a purely nonstandard one.
\end{remark}

Throughout, for each $\mathrm{e}\in E$, the constant function from $X$ into $E$ taking the value $\mathrm{e}$
will be denoted by $\tilde{\mathrm{e}}$.  Also, given a set $A$, $\chi_A$ stands 
for the characteristic function on $A$.

As usual, if there is no  confusion both the metric of $X$ and that of $Y$ will
be denoted by $d$.

Given a surjective linear isometry $T:\mathrm{Lip}(X,E)\rightarrow \mathrm{Lip}(Y,F)$, 
we denote 
$$\mathfrak{A} (T) :=\{y\in Y: T  \tilde{\mathrm{e}} (y)=0 \ \forall  \mathrm{e}\in E\}$$
and
$$\mathfrak{B} (T) := \mathfrak{A} (T)^c.$$
The partition of $Y$ into these two subsets will be very much used in Sections~\ref{purusa} and~\ref{asucar},
and the fact that $\at$ is empty will turn out to be basically equivalent to $T$ being standard. 
This property will receive a special name. 
 We define Property {\bf P} as follows:
\begin{description}\label{tak}
\item[P] For each $y \in Y$, there exists $\mathrm{e} \in E$ with $T  \tilde{\mathrm{e}} (y) \neq 0$.
\end{description}

\section{Main results}

 We first give some results ensuring that an isometry is standard, and then characterize spaces and describe
the isometries when this is not the case.   Theorem~\ref{marthahoyo} and Corollary~\ref{bocadillo} are proved in
 Section~\ref{rrollo}, and Theorem~\ref{nonelk} in Section~\ref{asucar}.

It is obvious that, by definition, if $T$ is {\em not} nonstandard, then it satisfies Property {\bf P}.
 The converse is given by Theorem~\ref{marthahoyo} and Corollary~\ref{bocadillo}.

\begin{theorem}\label{marthahoyo}
Let  $T:\mathrm{Lip}(X,E)\rightarrow \mathrm{Lip}(Y,F)$ be a surjective linear isometry satisfying Property {\bf P}. 
Then $E$ and $F$ are linearly isometric. 

Furthermore, if we are in any of the  following two cases:
\begin{enumerate}
\item\label{zala} $X$ and $Y$ are complete,
\item\label{mea} $E$ (or $F$) is not complete,
\end{enumerate}
then $T$ is standard.
\end{theorem}

\begin{remark}\label{subidahoyo}
In Theorem~\ref{marthahoyo}, we cannot in general ensure that the map $h$ 
is an isometry or that it preserves distances equal to $2$. 
Indeed,  following the same ideas as in \cite[Proposition 1.7.1]{W},
if $(Z, d)$ is a metric space with diameter  $\mathrm{diam}(Z, d) >2 $, then there is a new metric
$d' (\cdot, \cdot ) := \min \tl 2, d (\cdot, \cdot ) \tr$ on $Z$ with $\mathrm{diam}(Z, d') = 2$ such that $\mathrm{Lip}(Z, E)$ with respect to $d$
and $\mathrm{Lip}(Z, E)$ with respect to $d'$ are linearly isometric.
On the other hand, notice also that, if $d_1'$ and $d_2'$ are defined in a similar way, then the map 
$h: (Y, d_2) \ra  (X, d_1) $ belongs to $\iso$  if and only if
 $h: (Y, d_2') \ra    (X, d_1')$ is an isometry.
\end{remark} 

In Theorem~\ref{marthahoyo}, when (\ref{zala}) and (\ref{mea})  do not hold, $E$ and $F$ 
are complete and $X$ (or $Y$) is not. 
In this case, it is easy to see that in general $T$ is not standard. 
Nevertheless, 
 we have the following result.

\begin{corollary}\label{bocadillo}
Suppose that $E$ and $F$ 
are complete and $X$ or $Y$ is not. If $T:\mathrm{Lip}(X,E)\rightarrow \mathrm{Lip}(Y,F)$ is a surjective linear isometry satisfying Property {\bf P}, then $\widehat{T} :\mathrm{Lip} \pl \widehat{X},E \pr \rightarrow \mathrm{Lip} \pl \widehat{Y},F \pr $ is standard. 
\end{corollary}

 We next give the general form that a nonstandard isometry (or, equivalently, an isometry not satisfying Property {\bf P})
must take.

\begin{theorem}\label{nonelk}
Assume that we are in any of the  following two cases:
\begin{enumerate}
\item\label{kuris} $X$ and $Y$ are complete,
\item\label{masu} $E$ (or $F$) is not complete,
\end{enumerate}
Then there exists a nonstandard isometry $T:\mathrm{Lip}(X,E)\rightarrow \mathrm{Lip}(Y,F)$  
if and only if 
the following three conditions hold simultaneously:
\begin{enumerate}
\item $X$ and $Y$ are of type $\mathbf{A}$,
\item there exists $h \in \iso$,
\item $E$ and $F$ are linearly isometric.
\end{enumerate}

In this case, $T = S_{\varphi} \circ T'$, where $T' : \cxe \ra \cyf$ is a standard isometry   and 
$S_{\varphi} : \cyf \ra \cyf$ is a purely nonstandard isometry.
\end{theorem}

\begin{remark}\label{sinpase}
In the case when $E$ and $F$ 
are complete and $X$ or $Y$ is not,  if $T:\mathrm{Lip}(X,E)\rightarrow \mathrm{Lip}(Y,F)$ is a nonstandard isometry, then so is $\widehat{T}$, and the description given in Theorem~\ref{nonelk} applies to $\widehat{T}$. 
 In particular, $T = {\mathbf{i}_Y}^{-1} \circ  S_{\varphi} \circ T' \circ \mathbf{i}_X$
 where $T' : \cxeh \ra \cyfh$ is a standard isometry   and 
$S_{\varphi} : \cyfh \ra \cyfh$ is purely nonstandard. 
\end{remark}

A direct consequence (and easy to check) of Theorem~\ref{nonelk} and Remark~\ref{sinpase} is the following.

\begin{corollary}\label{deluxion}
If $E$ is not complete,  then there exists a nonstandard isometry from $\mathrm{Lip}(X,E)$ onto itself if and only if $X$ is of type $\mathbf{A}$. If $E$ is   complete,  then there exists a nonstandard isometry from $\mathrm{Lip}(X,E)$ onto itself if and only if $\widehat{X}$ is of type $\mathbf{A}$.
\end{corollary}

Theorem~\ref{nonelk} says that, under some assumptions, when two spaces of Lipschitz functions are linearly isometric, there
exists in fact a standard isometry between them. The following result is a simple consequence of Theorems~\ref{marthahoyo} and~\ref{nonelk}, Corollary~\ref{bocadillo} and Remark~\ref{sinpase}.

\begin{corollary}\label{rentreeflot}
$\mathrm{Lip}(X,E)$ and $\mathrm{Lip}(Y,F)$ are linearly isometric if and only if  $E$ and $F$ are linearly isometric
and 
\begin{itemize}
\item $\iso$ is nonempty (when $E$ and $F$ are not complete).
\item $\isol$ is nonempty (when $E$ and $F$ are   complete).
\end{itemize}  
\end{corollary}

We finally adapt the above results to the special case of metrics $d^{\alpha}$, $0 <\alpha <1$. Even if in this 
case we just deal with metrics and, consequently, the general form of the isometries between spaces $\cxea$ is 
completely given by Theorems~\ref{marthahoyo} and~\ref{nonelk}, 
Corollary~\ref{bocadillo} and Remark~\ref{sinpase}, it is interesting to see how  the condition of being of type $\mathbf{A}$
can be translated to metrics $d^{\alpha}$. This
 turns out to be more restrictive,  
and constitutes a generalization of
the scalar case on compact spaces given in \cite[Theorem 3.3]{M}.

\begin{definition}\label{rebotapel}
Let $0< \alpha <1$. We say that a metric space $(X,d)$ is of {\em type $\mathbf{A}_{\alpha}$} if there are a partition of $X$ into two subsets 
$\mathfrak{A}$, $\mathfrak{B}$, and a   map $\varphi : \mathfrak{A} \ra \mathfrak{B}$ with the following properties: 
\begin{enumerate}
\item $d(x, \varphi (x)) =1$ for every $x \in \ah$,
\item $d^{\alpha} (x, z) \ge 2$ whenever $x \in \mathfrak{A}$ and  $z \in \mathfrak{B}$, $z \neq \varphi (x)$,
\item $d^{\alpha} (x_1, x_2 ) \ge 2 $ whenever $x_1, x_2 \in \mathfrak{A}$ and $\varphi(x_1) \neq \varphi (x_2)$.
\end{enumerate}
\end{definition}

\begin{proposition}\label{pesafrank}
Let $0< \alpha <1$ and let $(X, d)$ be a metric space.
The following two statements are equivalent:
\begin{enumerate}
\item $\pl X, d^{\alpha} \pr$ is of type $\mathbf{A}$,
\item  $\pl X, d \pr$ is of type $\mathbf{A}_{\alpha}$.
\end{enumerate}
\end{proposition}

Proposition~\ref{pesafrank} will be proved in Section~\ref{asucar}.

It is clear that, since 
$X$ is of type $\mathbf{A}_{\alpha}$ if and only if its completion is, the statement of
Corollary~\ref{deluxion} is even simpler when dealing with $\cxea$.
 
Also, it is immediate to see that, if $(X, d)$ is of type $\mathbf{A}_{\alpha}$, then it is of type 
$\mathbf{A}_{\beta}$ for $\alpha < \beta \le 1$. Consequently, by Theorem~\ref{nonelk} and Remark~\ref{sinpase},
we conclude the following.

\begin{corollary}\label{subyecta}
Let $0< \alpha <1$. If there exists a nonstandard isometry between $\cxea$ and $\cyfa$, then there exists
a nonstandard isometry between $\mathrm{Lip}^{\beta} (X, E)$ and $\mathrm{Lip}^{\beta} (Y, F)$ whenever 
$\alpha < \beta \le 1$.
\end{corollary}

Obviously, the converse of Corollary~\ref{subyecta} is not true in general. The following example shows somehow the differences between  cases.
\begin{example}\label{retipnol}
Let $X := \{-1\} \cup (0,1) \subset \mathbb{R}$. $X$ is not of type $\mathbf{A}$, 
but its completion $\widehat{X} = \{-1\} \cup [0,1]$ is. Neither $X$ nor $\widehat{X}$ are of type 
$A_{\alpha}$ for $0 < \alpha <1$. Consequently, we have
\begin{itemize}
\item If $E$ is not complete, then all linear isometries from $\cxe$ onto itself are standard, but there are
nonstandard isometries from $\mathrm{Lip} \pl \widehat{X},E \pr$ onto itself.
\item If $E$ is  complete, then  there are nonstandard isometries from $\cxe$ onto itself. Obviously, by definition of nonstandard isometry, the same holds for 
$\mathrm{Lip} \pl \widehat{X},E \pr$. 
\item For every $E$ and $\alpha \in (0,1)$, all linear isometries from $\cxea$ onto itself are standard.
The same holds for $\mathrm{Lip}^{\alpha} \pl \widehat{X},E \pr$. In the case when $E$ is complete and $X$ is not, this is due to the special form of $X$.
\end{itemize}
\end{example}

\section{The case when $T$ satisfies Property {\bf P}}\label{rrollo}

In this section, unless otherwise stated, we assume that 
$T$ is a linear isometry from
$\mathrm{Lip}(X,E)$ onto $\mathrm{Lip}(Y,F)$ satisfying Property {\bf P}.

Our first goal consists of showing that $T$ is indeed an isometry with respect to the norm $\vc \cdot \vd_{\infty}$. The following two lemmas will be the key tools used to prove it.

\begin{lemma}\label{paco}
Let $f \in \mathrm{Lip}(X,E)$ and $x_0 \in X$ be such that $f(x_0) \neq 0$.
Then there exists $g \in \mathrm{Lip}(X,E)$ with
$\vc g(x_0) \vd = \vc g \vd_{\infty} > L(g)$
such that
$$\vc g (x_0) \vd + \vc f (x_0) \vd = \vc (g + f) (x_0) \vd = \vc g + f \vd_{\infty} > L (g + f).$$
\end{lemma}

\begin{proof}
We put $\mathrm{e}:=f(x_{0})$ and assume without loss of generality that $\vc \mathrm{e} \vd =1$.
We then consider $l \in \mathrm{Lip}(X,E)$ defined by
$
l (x):= \max \{0,2-d(x,x_{0})\} \cdot \mathrm{e}
$
for each $x\in X$. 
Clearly $l$ satisfies
$\vc l \vd_{L}=\vc l \vd_{\infty}=\vc l(x_{0}) \vd=2$ and $L(l)\leq 1$,
and also $\vc l(x) \vd< 2$ for all $x\in X$, $x\neq x_{0}$.

Take $n \in \mathbb{N}$ with $n > \vc f \vd_L$. Firstly, it is easy to check that
$\vc (nl +f)(x_{0}) \vd=2n+ 1$. On the other hand, we also have that
$L(nl+f)\leq nL(l)+L(f)< 2n$ and, for $x\in X \setminus \tl x_{0} \tr$ with $d(x, x_0) <2$,
\begin{eqnarray*}
\vc (nl+f)(x) \vd &\leq& n\vc l(x) \vd + \vc f(x) \vd \\
&\leq& 2n-nd(x,x_{0})+\vc f(x_{0}) \vd +  L(f) d(x,x_{0}) \\
&<& 2n+1,
\end{eqnarray*}
whereas if $d(x, x_0 ) \ge 2$, then $l(x)=0$ and
$\vc (nl+f)(x) \vd \leq  \vc f \vd_L < n$.

Consequently, if we define $g:= nl$, the lemma is proved.
\end{proof}

\begin{remark}\label{guelga}
It is easy to see that if $f_1 , f_2 \in \mathrm{Lip}(X,E)$ satisfy $\vc f_1 (x_0) \vd , \vc f_2 (x_0) \vd <  \vc f (x_0) \vd$, then the proof of Lemma~\ref{paco} can be slightly modified (by taking 
$n >  \vc f \vd_L, \vc f_1 \vd_L , \vc f_2 \vd_L $) so   that 
$\vc g + f_i \vd_{\infty} < \vc g + f \vd_{\infty}$ for $i=1, 2$.
\end{remark}

\begin{lemma}\label{1750}
If $f \in \mathrm{Lip}(X,E)$ satisfies $\vc Tf (y_0) \vd = \vc Tf \vd_{\infty} > L( Tf )$
for some $y_0 \in Y$, then $L(f) \le \vc f \vd_{\infty}$.
\end{lemma}

\begin{proof}
Suppose that $\vc  f \vd_{\infty} < L(f)$. Then,
for each $\mathrm{e} \in E$, there exists $M>0$ such that 
$\vc f \vd_{\infty} + M \vc \mathrm{e} \vd < L(f) $,  
so $\vc f \pm M \tilde{\mathrm{e}} \vd_{\infty}
< L(f) = L(f \pm M \tilde{\mathrm{e}})$.
 Therefore,
$$
\vc f \pm M \tilde{\mathrm{e}} \vd_{L}=
L(f \pm M \tilde{\mathrm{e}})=L(f)=\vc f \vd_{L}.
$$
Since $T$ is an isometry,
$\vc Tf \pm M T\tilde{\mathrm{e}} \vd_{L}=\vc Tf \vd_{L} = \vc Tf (y_0)  \vd$,
which implies in particular that
$\vc Tf (y_0)  \pm M T\tilde{\mathrm{e}}(y_{0}) \vd\leq\vc Tf (y_{0}) \vd$ and, by Inequality~(\ref{strictly}), 
 that $T \tilde{\mathrm{e}} (y_0) =0$ for every $\mathrm{e} \in E$,
which  goes against our hypotheses.
\end{proof}

\begin{remark}\label{hirche}
Notice that, in the proof of Lemma~\ref{1750}, we just use the fact that there exists $\mathrm{e} \in E$ with $T  \tilde{\mathrm{e}} (y_0) \neq 0$, and not the general assumption that Property {\bf P} holds.
\end{remark}

\begin{corollary}\label{mariusz}
$T$ is an isometry with respect to the supremum norm.
\end{corollary}

\begin{proof}
Assume that $\vc f \vd_{\infty} < \vc Tf \vd_{\infty}$, and pick $\epsilon >0$ and $y_0 \in Y$
such that $\vc f \vd_{\infty} + \epsilon < \vc Tf (y_0) \vd$. Next, by Lemma~\ref{paco}, we can take
$g \in \mathrm{Lip}(Y,F)$ with $\vc g(y_0) \vd = \vc g \vd_{\infty} > L(g)$ and such that
$$\vc g (y_0) \vd + \vc Tf (y_0) \vd = \vc (g + Tf) (y_0) \vd = \vc g + Tf \vd_{\infty} > L (g + Tf).$$
Applying Lemma \ref{1750}, we conclude both that 
$$L(T^{-1} g) \le \vc  T^{-1} g \vd_{\infty} = \vc g (y_0) \vd$$ and
 that $$L(T^{-1} g+f) \le \vc T^{-1} g+f  \vd_{\infty} = \vc g + Tf \vd_L = \vc (g + Tf) (y_0) \vd .$$
  But this is impossible because
$$
\vc T^{-1} g+ f  \vd_{\infty} < \vc g (y_0) \vd+\vc Tf (y_0) \vd - \epsilon < \vc (g + Tf) (y_0) \vd.
$$
We conclude that $\vc Tf \vd_{\infty} \le \vc f \vd_{\infty}$ for every $f$. 
We next see that $T^{-1}$ also satisfies Property {\bf P}, which is enough to prove the equality. By the above, given a nonzero $\mathrm{f} \in F$   
we have  $\vc T^{-1} \tilde{\mathrm{f}} \vd_{\infty} = \vc \mathrm{f} \vd$. Also,
 if $\pl  T^{-1} \tilde{\mathrm{f}} \pr (x_0) =0$ for some $x_0 \in X$, then there exists  
$k \in \mathrm{Lip}(X,E)$, $k \neq 0$, with $\vc k (x) \vd + \vc  \pl T^{-1} \tilde{\mathrm{f}} \pr (x) \vd \le \vc \mathrm{f} \vd$
for every $x \in X$. By Inequality~(\ref{strictly}), $\vc \mathrm{f} \vd < \vc \tilde{\mathrm{f}} + T k \vd_{\infty}$
or $\vc \mathrm{f} \vd < \vc \tilde{\mathrm{f}} - T k \vd_{\infty}$, which contradicts the paragraph above.
\end{proof}

\begin{remark}\label{hadela}
Notice that, in the proof of Corollary~\ref{mariusz}, we have  seen that  $T^{-1}$ also satisfies Property {\bf P}.
\end{remark}

We are   now ready to see that, under the assumptions we make in this section, every surjective 
linear isometry  
 is biseparating.   
 
\begin{proposition}\label{bis}
$T $ is biseparating.
\end{proposition}

\begin{proof}
We prove that $T$ is separating. Suppose that it is not, so 
there exist $f,g\in \mathrm{Lip}(X,E)$
such that $c(f)\cap c(g)=\emptyset$ but $Tf(y_{0})=\mathrm{f}_1\neq 0$ and $Tg(y_{0})=\mathrm{f}_2\neq 0$
for some $y_{0}\in Y$. Taking into account Inequality~(\ref{strictly}), we can assume without loss of
generality that $\vc \mathrm{f}_2 \vd \leq  \vc \mathrm{f}_1 \vd < \vc \mathrm{f}_1+ \mathrm{f}_2 \vd$.
Now, 
by Lemma~\ref{paco} and Remark~\ref{guelga}, there exists $k \in \mathrm{Lip}(Y,F)$ such that 
$\vc k + Tf \vd_{\infty} , \vc k + Tg \vd_{\infty} < \vc k + Tf + Tg \vd_{\infty}$.

On the other hand,  
since $f$ and $g$
have disjoint cozeros,  
$$
\vc T^{-1} k+ f+ g \vd_{\infty}=\max \tl \vc T^{-1} k+ f \vd_{\infty},\vc T^{-1} k+  g \vd_{\infty} \tr,
$$
and consequently 
$
\vc k+Tf+Tg \vd_{\infty}=\max\tl \vc k+Tf \vd_{\infty},\vc k+Tg \vd_{\infty} \tr
$,
which is a contradiction.  

By Remark~\ref{hadela}, $T^{-1}$ is also separating.
\end{proof}

\begin{remark}\label{pinbo}
In \cite[Theorem 3.1]{AD} (see also comments after it)  a description of biseparating maps
$S: \mathrm{Lip}(X,E)\rightarrow \mathrm{Lip}(Y,F)$ is given, but we cannot use it here because 
assumptions of completeness on $X$ and $Y$ are made in \cite{AD}. Under some circumstances automatic continuity
of such $S$ can be achieved and, in that case, the description goes  as follows 
 (where $\mathcal{L}(E,F)$ denotes 
 the normed space of all
linear and continuous operators from $E$ to $F$):
There exist a  homeomorphism $k: Y \rightarrow X$ 
and  a   map   $K : Y \rightarrow \mathcal{L}(E, F)$ 
(which is easily seen to be  also Lipschitz with 
$L(K) \le \vc S \vd$) such that
$
Sf(y)=Ky (f(k(y)))
$
for all $f\in \mathrm{Lip}(X,E)$ and $y\in Y$. Also, if both $X$ and $Y$ are bounded, then the map $k$
is bi-Lipschitz.
\end{remark}

\begin{proposition}\label{nuevaescocia}
Given $ \mathrm{e}\in E$, $T \tilde{\mathrm{e}}$ is constant on each $1$-component of $Y$ and 
$\vc T \tilde{\mathrm{e}} (y) \vd = \vc \mathrm{e} \vd$ for all $y \in Y$.
\end{proposition}

\begin{proof}
Suppose that this is not the case, but there exist $\mathrm{e}  \in E$, $\vc \mathrm{e}  \vd =1$,
and $y_1, y_2 \in Y$ with $y_1 \sim_1 y_2$, such that $\mathrm{f}_1 := T \tilde{\mathrm{e}} (y_1)$
and $\mathrm{f}_2 := T \tilde{\mathrm{e}} (y_2)$ are different. Of course,
 we may assume without loss of generality that $D:= d(y_1, y_2) <1$
 and that 
$ \mathrm{f}_1  \neq 0$. Now, if we  consider $g \in \mathrm{Lip}(Y,F)$ defined by
$$
g ( y):= \max\tl 0,1 - \frac{d(y, y_1 )}{D} \tr \cdot \mathrm{f}_1
$$
for all $y\in Y$, then obviously 
$\vc g \vd_L = L (g) =  \vc \mathrm{f}_1 \vd / D > \vc g \vd_{\infty}$.
As a consequence, using Corollary~\ref{mariusz}, $\vc T^{-1} g \vd_L > \vc T^{-1} g \vd_{\infty}$,
and we can take $M>0$ such that
$\vc   T^{-1} g \pm M \tilde{\mathrm{e}} \vd_L = \vc T^{-1} g \vd_L = \vc \mathrm{f}_1 \vd /D$.

Notice also that, as $F$ is strictly convex, either
$\vc \mathrm{f}_1 + M \pl \mathrm{f}_1 - \mathrm{f}_2 \pr  \vd  > \vc \mathrm{f}_1 \vd $
or
$\vc \mathrm{f}_1 - M \pl \mathrm{f}_1 - \mathrm{f}_2 \pr  \vd  > \vc \mathrm{f}_1 \vd $,
that is, either
$$\frac{\vc \pl g + M T \tilde{\mathrm{e}} \pr   (y_1) -
\pl g + M T \tilde{\mathrm{e}} \pr   (y_2)\vd }{d(y_1, y_2)} > \frac{\vc \mathrm{f}_1 \vd}{D}$$
or
$$\frac{\vc   \pl g - M T \tilde{\mathrm{e}} \pr   (y_1) -
\pl g - M T \tilde{\mathrm{e}} \pr   (y_2)\vd }{d(y_1, y_2)} > \frac{\vc \mathrm{f}_1 \vd}{D},$$
which implies that either
$\vc g + M T \tilde{\mathrm{e}}  \vd_L > \vc \mathrm{f}_1 \vd /D$
or
$\vc g - M T \tilde{\mathrm{e}} \vd_L > \vc \mathrm{f}_1 \vd /D$,
yielding a contradiction.

Finally suppose that $T \tilde{\mathrm{e}} (y) = \mathrm{f}$ for all $y $ in  $B \in \mathtt{Comp}_1 (Y)$. 
Since  by Proposition~\ref{bis} $T^{-1}$ is separating,   
$c \pl T^{-1} \pl \chi_B \cdot T \tilde{\mathrm{e}} \pr \pr \cap c \pl T^{-1} \pl \chi_{Y \setminus B} \cdot T \tilde{\mathrm{e}} \pr \pr = \emptyset$. This implies that $\mathrm{e}$ is   
the only nonzero  value taken by $T^{-1} \pl \chi_B \cdot T \tilde{\mathrm{e}} \pr$ and, since $T^{-1}$ is an isometry with respect to $\vc \cdot \vd_{\infty}$, we have that $\vc \mathrm{e} \vd = \vc \mathrm{f} \vd$.
\end{proof}

\begin{lemma}\label{lavaderomartha}
There exists a bijection $H : \mathtt{Comp}_1 (X) \ra \mathtt{Comp}_1 (Y)$ and, for each $A \in \mathtt{Comp}_1 (X)$,
a surjective linear isometry $\mathbf{J}_A : E \ra F$ with the property
that  
$T \pl \chi_A \cdot \tilde{\mathrm{e}} \pr = \chi_{H(A)} \cdot \widetilde{\mathbf{J}_A (\mathrm{e})} = \chi_{H(A)} \cdot T    \tilde{\mathrm{e}}   $ for every 
$\mathrm{e} \in E$.
\end{lemma}

\begin{proof}
Fix $A \in \mathtt{Comp}_1 (X)$ and $\mathrm{e} \in E$ with $\vc \mathrm{e} \vd =1$, and take
 $f := \chi_A \cdot \tilde{\mathrm{e}}$, $g := \chi_{X \setminus A} \cdot \tilde{\mathrm{e}}$ in $\mathrm{Lip}(X,E)$. 
We have that $c(f) \cap c(g) = \emptyset$,  so by Proposition~\ref{bis} $Tf$ and $Tg$ have disjoint cozeros. 
Then, by Proposition~\ref{nuevaescocia}, $Tf(y), Tg(y) \in \tl 0, T \tilde{\mathrm{e}} (y) \tr$
for all $y  \in Y$.
Now, suppose that $y \sim_1 y'$,    
and that $Tf(y) \neq 0$ and
$Tf (y') =0$. We
 can assume that $d(y, y') <1$. 
Since $\vc Tf (y) \vd =1$, we deduce $ L(Tf) \ge 1/ d(y,y') >1$, which
is impossible because $\vc f \vd_L =1$.

Reasoning similarly with $T^{-1}$,  
$Tf = \chi_B \cdot \tilde{\mathrm{f}}$ 
for some $1$-component $B$ in $Y$ and some norm-one vector $\mathrm{f} \in F$. The conclusion is now easy.
\end{proof}

\begin{lemma}\label{stomas}
Given $A, B \in \mathtt{Comp}_1 (X)$, if $\min \tl d(A, B) , d(H(A), H(B)) \tr <2$, then 
$d(A, B) = d(H(A), H(B))$ and $\mathbf{J}_A   = \mathbf{J}_B$.
\end{lemma}

\begin{proof}
Put $D_1 := d  (A, B)$, $ D_2 := d  (H(A), H(B))$. 
 Due to the symmetric r\^{o}les
of $H$ and $H^{-1}$  with respect to $T$ and $T^{-1}$, we can assume without loss of generality that 
$D_1 \le D_2$.
Pick $\mathrm{e} \in E$ with $\vc \mathrm{e} \vd =1$, and define
 $ f := \pl \chi_A  - \chi_B \pr \cdot \tilde{\mathrm{e}} \in \mathrm{Lip}(X,E)$. We easily see that
 $\vc f \vd_L = L(f) = 2/D_1$  and, since 
$L (Tf) = \vc \mathbf{J}_A (\mathrm{e} ) + \mathbf{J}_B (\mathrm{e} ) \vd  / D_2 
\le 2/D_2$,  
we necessarily have    $D_1 = D_2$ and  
$\vc \mathbf{J}_A (\mathrm{e} ) + \mathbf{J}_B (\mathrm{e} ) \vd =2$, so $\mathbf{J}_A (\mathrm{e}  ) = \mathbf{J}_B (\mathrm{e} )$ because $F$ is strictly convex. 
\end{proof}

\begin{corollary}\label{borjas}
There exists a map $J: Y \ra I(E, F)$ which is constant on each $2$-component of $Y$ and such that
$T \tilde{\mathrm{e}} (y) = J y (\mathrm{e})$ for all $\mathrm{e} \in E$ and $y \in Y$.
\end{corollary}

\begin{proof}
We define $Jy:= \mathbf{J}_A $ if $y \in H(A)$ and $A \in \mathtt{Comp}_1 (X)$. 
Applying Lemma~\ref{stomas}, 
 the result follows.
\end{proof}

\begin{lemma}\label{gatosinchicle}
Let $(y_n)$ be a  sequence in $Y$ which is not a Cauchy sequence and such that
all $y_n$ are pairwise different. 
Then there exist infinite   subsets $A_1$ and $A_2$
of $\tl y_n : n \in \mathbb{N} \tr$ with  $d \pl A_1, A_2 \pr >0$.
\end{lemma}

\begin{proof}
Taking a subsequence if 
necessary, 
we have that there exists $\epsilon>0$ such that
 $d (y_{2n} , y_{2n-1} ) \ge 3 \epsilon $ for all 
$n  \in \mathbb{N}$. Let $  A:= \tl y_n : n \in \mathbb{N} \tr$. Now we have two possibilities: either there exists $n_0$ such that 
$B \pl y_{n_0} ,  \epsilon \pr$ contains infinitely many $y_n$ or $A \cap B \pl y_k ,  \epsilon \pr$ is
 finite for every $k$. In the first case, it is clear that $A_1 := A \cap B \pl y_{n_0} ,  \epsilon \pr$
 and $A_2 := \tl y_{2n} : y_{2n -1} \in A_1 \tr \cup \tl y_{2n-1} : y_{2n} \in A_1 \tr$ satisfy 
$d(A_1, A_2) \ge \epsilon$. In the second case,   we can find a subsequence $\pl y_{n_k} \pr$
with $d \pl y_{n_k}, y_{n_l} \pr > \epsilon$ when $k \neq l$, and the result follows easily. 
\end{proof}

 In Lemma~\ref{uberrimo} and Corollary~\ref{estela} we do not necessarily assume that base spaces are not complete, so it could be the case that $\widehat{X} =X$ and $\widehat{Y} =Y$. 

\begin{lemma}\label{uberrimo}
Given 
$x_0 \in \widehat{X}$, there exists  
$y_0 \in \widehat{Y}$ such that $\widehat{ Tf } (y_0) = 0$ whenever $f \in \mathrm{Lip}(X,E)$ satisfies 
$\widehat{ f } (x_0) =0$. 
\end{lemma}

\begin{proof}
Fix $\mathrm{e} \in E$ with $\vc \mathrm{e} \vd =1$, and let 
$$A := \tl f \cdot \tilde{\mathrm{e}} : f \in \mathrm{Lip}(\widehat{X}) , \ f (x_0) =1 , \ \forall \epsilon > 0 \sup_{d(x , x_0) \ge \epsilon} \va f(x) \vb <1   \tr .$$
We will see that there exists a unique point $y_0 \in \widehat{Y}$ such that 
 $\widehat{\vc Tf \vd} (y_0) =1$ for every $f \in A$. 

Fix $f_0 \in A$. By Corollary~\ref{mariusz}, taking into account that $\vc f_0 \vd_{\infty} =1$, 
there exists a sequence $(y_n)$ in $Y$ such that $\vc Tf_0 (y_n) \vd \ge 1 - 1/n$ for each $n \in \mathbb{N}$. 
Let us see that it is a Cauchy sequence. Suppose that this is not the case. Either if all $y_n$ are pairwise different (by using Lemma~\ref{gatosinchicle})
 or not, we see that  there exist   subsets $A_1, A_2$ of $\tl y_n : n \in \mathbb{N} \tr$ such that
$d(A_1, A_2) >0$ and  $\sup_{y_n \in A_i} \vc Tf_0 (y_n) \vd =1$, $i=1,2$. Then we take
$g_1, g_2 \in \mathrm{Lip}(\widehat{Y})$ with  $0 \le g_1, g_2 \le 1$ such that
$g_1 (A_1) = 1$, $g_2  (A_2) = 1$, and $g_1 g_2 \equiv 0$. It is immediate that $\vc Tf_0 + g_i Tf_0 \vd_{\infty}  =2$ 
for $i=1,2$. 
 Since, again by Corollary~\ref{mariusz},
 $\vc T^{-1} \pl g_i Tf_0 \pr \vd_{\infty} =1$,
 we 
deduce that $ \widehat{\vc  T^{-1} \pl g_i Tf_0 \pr  \vd} (x_0) =1$ for $i=1,2$, which goes against the fact that
$T^{-1}$ is separating. Consequently $(y_n)$ is a Cauchy sequence and converges to a point $y_0 \in \widehat{Y}$,
which obviously satisfies $\widehat{\vc Tf_0 \vd} (y_0) =1$. Now it is straightforward to see that
$\widehat{\vc Tf \vd} (y_0) =1$ for every $f \in A$.

Next suppose that $f \in \mathrm{Lip}(X, E)$ satisfies
 $\widehat{f} (x_0)=0$. Then, given $\epsilon >0$, there exists 
$f_{\epsilon} \in \mathrm{Lip}(X, E)$ such that $\widehat{f_{\epsilon}} \equiv 0$
 on a neighborhood of $x_0$ and 
$\vc f - f_{\epsilon} \vd_{\infty} \le \epsilon$. We can  take $f_{\epsilon}' \in A$ with
 $c(f_{\epsilon}') \cap c(f_{\epsilon}) = \emptyset$, and we
deduce from the paragraph above that  $\widehat{\vc Tf_{\epsilon} \vd} \equiv 0$ on a neighborhood of $y_0$ in $\widehat{Y}$; in 
particular $\widehat{\vc Tf_{\epsilon}  \vd }  (y_0) =0$. Since $\vc T f -  T f_{\epsilon} \vd_{\infty} \le \epsilon$ (by Corollary~\ref{mariusz}),
 we conclude that $ \widehat{ \vc Tf   \vd} (y_0) \le \epsilon$, and we are done.
\end{proof}

\begin{corollary}\label{estela}
There exists a bijective map $h : \widehat{Y} \ra \widehat{X}$ such that
 $ Tf   (y) = J y (\widehat{f } (h (y)))$  
whenever  $y \in Y$ 
and $f \in \mathrm{Lip}(X,E)$ admits a continuous extension to $h(y)$.
\end{corollary}

\begin{proof}
Let $x_0$ and $y_0$ be as in Lemma~\ref{uberrimo}. Since $T^{-1}$ is also biseparating,  
there exists $x_1 \in \widehat{X}$ such that $ \widehat{ f } (x_1) =0$ whenever $\widehat{ Tf } (y_0) = 0$ and,
in particular, whenever $ \widehat{ f } (x_0) =0$. Now, as $\mathrm{Lip}(\widehat{X}, E)$   separates points
in $\widehat{X}$, we deduce that $x_1 = x_0$. 
As a consequence, it is straightforward to see that Lemma~\ref{uberrimo} gives us a bijective map between
 $\widehat{X}$ and $\widehat{Y}$, which we denote by $h: \widehat{Y} \ra \widehat{X}$, satisfying 
 $\widehat{ Tf  } (y) =0$  if and only if  $\widehat{ f  } (h(y)) =0$.
Finally,
if $f \in \mathrm{Lip}(X,E)$ can be continuously extended to $h(y)$, say
$\widehat{f } (h(y)) = \mathrm{e} \in E$, then
$\widehat{\pl f - \tilde{\mathrm{e}}  \pr } (h(y)) =0$,
and the representation follows from  
Corollary~\ref{borjas}.
\end{proof}

\begin{remark}\label{pliniuz}
As in the proof of Corollary~\ref{estela}, the bijection $k: \widehat{X} \ra \widehat{Y}$ associated to $T^{-1}$ satisfies $\widehat{ T^{-1} g  } (x) =0$  if and only if  $\widehat{ g  } (k(x)) =0$,
$g \in \mathrm{Lip}(Y,F)$. This implies that $k = h^{-1}$.

\end{remark}

\begin{lemma}\label{ester}
If $E$ is not complete, then there exists a sequence $(\mathrm{e}_n)$ in $E$
with $\vc \mathrm{e}_n \vd \le 1/4^n$ such that $\sum_{n=1}^{\infty} \mathrm{e}_n$ does not converge in $E$.
\end{lemma}

\begin{proof}
Clearly, there exists a nonconvergent sequence $(\mathrm{u}_n)$ in $E$ satisfying 
$\vc \mathrm{u}_n - \mathrm{u}_{n+1} \vd \le 1/4^n$ for every
 $n \in \mathbb{N}$. It is then easy to check that it is enough to define
$ \mathrm{e}_n := \mathrm{u}_n - \mathrm{u}_{n+1}$ for each $n$.
\end{proof}

\begin{corollary}\label{sara}
If $E$  is not complete,
then the map $h$ given in
Corollary~\ref{estela} is a bijection from $Y$ onto $X$.
\end{corollary}

\begin{proof}
We will prove  first that $h(y) \in X$ whenever $y \in Y$. If this is not the case,
then take $y \in Y$ with $h(y) \in \widehat{X} \setminus X$. For each $n \in \mathbb{N}$, let
$$
f_n (x) := \max \tl 0,1 - 2^n d \pl x,   h(y) \pr   \tr 
$$
for all $x \in   X$. It is clear that each $f_n$ belongs to  
$\mathrm{Lip}(X)$
and that $L(f_n) \le  2^n$. It is easy to see that, since $\mathrm{Lip}(X,\widehat{E})$ is complete, if we take $(\mathrm{e}_n)$ in $E$ as in Lemma~\ref{ester}, then 
$f:= \sum_{n=1}^{\infty} f_n \cdot \widetilde{\mathrm{e}_n}$ belongs to 
$\mathrm{Lip}(X,\widehat{E})$, and since all values are taken in $E$, to  $\mathrm{Lip}(X,E)$.
Thus, by Corollary~\ref{mariusz},
$$
\lim_{k \ra \infty} \vc Tf - \sum_{n=1}^k T \pl f_n \cdot  \widetilde{\mathrm{e}_n} \pr \vd_{\infty}=0.
$$
Finally, by Corollary~\ref{estela}, this implies that
$Tf (y) = \sum_{n=1}^{\infty} J y (\mathrm{e}_n)$, which belongs to $\widehat{F} \setminus  F$,
and $Tf$ takes values outside $F$, which is absurd.

We  deduce  from Remark~\ref{pliniuz} that $h(Y) =X$.
\end{proof}

\begin{proof}[Proof of Theorem~\ref{marthahoyo}]
Taking into account Corollaries~\ref{borjas}, \ref{estela} and \ref{sara}, it is enough  to show that $h \in \iso$. Let $y_1, y_2 \in Y $ be such that $d(y_1, y_2) < 2$. 
We are going to see that $D:= d(h(y_1), h(y_2)) \le d(y_1, y_2)$.

Pick $\mathrm{e} \in E$ with $\vc \mathrm{e} \vd=1$ and define $g \in \mathrm{Lip}(X,E)$ by $$
g ( x):= \max\tl -1 ,1 - \frac{ 2 \ d(x, h(y_1) )}{D} \tr \cdot \mathrm{e}
$$
for every $x \in X$. We have that $\vc g \vd_{\infty}=1$, $L(g) = 2/D$, $ g (h(y_1)) = \mathrm{e} $,
and $ g (h(y_2)) = - \mathrm{e}$. Obviously, by Corollary~\ref{borjas}, $J y_1 = J y_2$, and 
$$
1 < \frac{2}{d(y_1, y_2)} =  \frac{\vc J y_1 (\mathrm{e})  - J y_2 (- \mathrm{e}) \vd}{d(y_1, y_2)}
 =   \frac{\vc T g (y_1) - T g (y_2) \vd}{d(y_1, y_2)} \le \vc Tg \vd_L,
$$
which implies that $\vc g \vd_L >1$, and then  $ \vc g \vd_L = L(g) = 2/D $. This means that
$\vc Tg \vd_L =2/D$, and consequently
$
2/d(y_1, y_2)  \le  2/D
$.
The other inequality can be seen in a similar way working with $T^{-1}$ (see Remark~\ref{pliniuz}).
\end{proof}

\begin{proof}[Proof of Corollary~\ref{bocadillo}]
The fact that $\widehat{T}$ satisfies Property {\bf P} follows easily from Proposition~\ref{nuevaescocia}. The conclusion is then immediate by Theorem~\ref{marthahoyo}. 
\end{proof}

\section{The distance between $\at$ and $\bt$}\label{purusa}

Propositions~\ref{siegediesref} and~\ref{floater} will be used in Section~\ref{asucar}.

\begin{proposition}\label{siegediesref}
Let $T:\mathrm{Lip}(X,E)\rightarrow \mathrm{Lip}(Y,F)$ be a surjective linear isometry. If $ \mathfrak{A} (T) \neq \emptyset$, then $d \pl \mathfrak{A} (T) , \mathfrak{B} (T) \pr \ge 1$.
\end{proposition}

\begin{proof}
Obviously $ \mathfrak{B} (T) \neq \emptyset$. Suppose first that $d  ( \mathfrak{A} (T),  \mathfrak{B} (T))<1$, and  take
$y_{0}\in  \mathfrak{B} (T)$ and $\epsilon >0$ with
$d  (y_{0},  \mathfrak{A} (T))<1-2\epsilon$.
We then  select $\mathrm{f}\in F$, $\vc \mathrm{f} \vd=1$, and define  $ l \in \mathrm{Lip}(Y,F)$ by
$
l(y):= \max\{0,2-d(y,y_{0})\}\cdot \mathrm{f}
$
for every $y\in Y$. We have that  
$\vc l \vd_{L}=\vc l \vd_{\infty}=\vc l(y_{0}) \vd=2$, $L(l)\leq 1$,
and  $\vc l(y) \vd< 2$ for all $y\in Y \setminus \tl y_0 \tr$.

Now, by Lemma~\ref{1750} (see also Remark~\ref{hirche}), we have that $L(T^{-1}l) \le \vc T^{-1}l \vd_{\infty}$.
Consequently $\vc T^{-1}l \vd_{L}=\vc T^{-1}l \vd_{\infty}$, and then
$\vc T^{-1}l \vd_{\infty}=2$. Therefore,  there is a point $x_{0}$ in $X$ such that
$
\vc T^{-1}l(x_{0}) \vd> 2-\epsilon
$,
that is, $T^{-1}l(x_{0})=\alpha \mathrm{e}$ for some  $\mathrm{e}\in E$, $\vc \mathrm{e} \vd=1$,
and $\alpha\in \mathbb{R}$, $\alpha>2-\epsilon$. Next, obviously
$$
\vc \tilde{\mathrm{e}}+T^{-1}l \vd_{L}
\geq \vc \mathrm{e}+T^{-1}l(x_{0}) \vd
= \vc \pl 1 + \alpha  \pr \mathrm{e} \vd 
 > 3- \epsilon,
$$
so $\vc T\tilde{\mathrm{e}}+ l \vd_{L} > 3-\epsilon$.
Since $L(T\tilde{\mathrm{e}}+ l)\leq L(T\tilde{\mathrm{e}})+L(l)\leq 2$,
this implies that $\vc T\tilde{\mathrm{e}}+ l \vd_{\infty} > 3-\epsilon$, and hence 
  the set 
$B := \tl y \in Y :  \vc \pl T \tilde{\mathrm{e}}+ l \pr (y) \vd > 3-\epsilon \tr$ is nonempty.

Notice that, since
$\vc T \tilde{\mathrm{e}} \vd_{\infty}  \le 1$,
all points $y \in B$ must satisty $\vc l(y) \vd >2-\epsilon$, which is equivalent to $d(y,y_{0})<\epsilon$.
Thus, for some $y_{1}$ with $d(y_{1}, y_{0})  < \epsilon$, 
we have $\vc T\tilde{\mathrm{e}}(y_{1})+ l(y_{1}) \vd > 3-\epsilon$, which implies that
$\vc T\tilde{\mathrm{e}}(y_{1})\vd > 1-\epsilon$.
On the other hand, taking into account that $d  (y_{0}, \mathfrak{A} (T) )<1-2\epsilon$, there exists $y_{2}\in \mathfrak{A} (T) $ with
$d(y_{0},y_{2})\leq 1-2\epsilon$. Finally, observe that
$$
\frac{\vc T\tilde{\mathrm{e}}(y_{1})-T\tilde{\mathrm{e}}(y_{2})\vd}{d(y_{1},y_{2})}
=\frac{\vc T\tilde{\mathrm{e}}(y_{1})\vd}{d(y_{1},y_{2})}
>\frac{1-\epsilon}{1-2\epsilon+\epsilon}=1,
$$
which allows us to conclude that $L(T\tilde{\mathrm{e}})>1$, in contradiction with
the fact that $\vc \mathrm{e} \vd=1$ and $T$ is an isometry.
\end{proof}

\begin{proposition}\label{floater}
Let $T:\mathrm{Lip}(X,E)\rightarrow \mathrm{Lip}(Y,F)$ be a surjective linear isometry. If $ y_0 \in \mathfrak{A} (T) $, then $d \pl y_0 , \mathfrak{B} (T) \pr = 1$.
\end{proposition}

\begin{proof}
Suppose on the contrary  that  there exists
$s \in (0, 1)$ such that $d \pl B (y_0, s) , \mathfrak{B} (T) \pr > 1+s$. Take 
 $f \in \mathrm{Lip}(Y)$ with $c(f) \subset B(y_0, s)$ and such that $0 \le f \le s$, $f(y_0) =s$, and 
$L (f) \le 1$. Let $\mathrm{e} \in E$ and $\mathrm{f} \in F$ have norm $1$.  
It is easy to check that
$\vc f \cdot  \mathrm{f} \pm   T\tilde{\mathrm{e}} \vd \le 1$, whereas, since 
$T^{-1} \pl f \cdot  \mathrm{f} \pr \neq 0$, Inequality~(\ref{strictly}) implies that 
$$\vc T^{-1} \pl f \cdot  \mathrm{f} \pr + \tilde{\mathrm{e}} \vd_{\infty} >1$$ or 
$$\vc T^{-1} \pl f \cdot  \mathrm{f} \pr - \tilde{\mathrm{e}} \vd_{\infty} >1 ,$$
contradicting the fact that $T$ is an isometry.
\end{proof}

We next see that Property {\bf P} 
holds when $Y$  
is $1$-connected.  Obviously, the same result holds if $X$ is $1$-connected (see Remark~\ref{hadela}).

\begin{corollary}\label{siegediesrefe}
 Let $Y$  be $1$-connected and suppose that $\mathrm{Lip}(X,E) $ and 
$ \mathrm{Lip}(Y,F)$ are linearly isometric. Then $X$ is also $1$-connected and
every surjective linear isometry $T:\mathrm{Lip}(X,E)\rightarrow \mathrm{Lip}(Y,F)$ 
satisfies Property {\bf P}.
\end{corollary}

\begin{proof}
By Proposition~\ref{siegediesref}, Property {\bf P} 
holds when $Y$  
is $1$-connected.  
 
The fact that $X$ is $1$-connected can be easily deduced from the representation of $T$ in
Theorem~\ref{marthahoyo} or that of $\widehat{T}$ in Corollary~\ref{bocadillo} (taking into account that
a metric space is $1$-connected if and only if so is its completion).
\end{proof}

\begin{remark}\label{guelgafunc}
An immediate consequence of Corollary~\ref{siegediesrefe} is that, when $X$ (or $Y$) is $1$-connected, every surjective linear isometry $T:\mathrm{Lip}(X,E)\rightarrow \mathrm{Lip}(Y,F)$
is standard in any of the cases (\ref{zala}), (\ref{mea}),  given in Theorem~\ref{marthahoyo}.
\end{remark}

\section{The case when $T$ does not satisfy Property {\bf P}}\label{asucar}

In this section, unless otherwise stated, we assume that 
$T$ is a linear isometry from
$\mathrm{Lip}(X,E)$ onto $\mathrm{Lip}(Y,F)$ that {\em does not} satisfy Property {\bf P} (that is, $\mathfrak{A} (T)  \neq \emptyset$). We will make use of Theorem~\ref{marthahoyo}, so we  also assume  that
 we are in any of the  following two cases:
\begin{enumerate}
\item $X$ and $Y$ are complete,
\item $E$ (or $F$) is not complete.
\end{enumerate}
It is then clear by Proposition~\ref{siegediesref} that $X$ is complete 
if and only if both $\as$ and $\bs$ are complete.

We will introduce two isometries on spaces of Lipschitz functions
defined on $\as$ and $\bs$. The fact that these new isometries turn out to be standard
  will allow us to obtain a description of $T$.

\begin{lemma}\label{1750ref}
Suppose  that  $f \in \mathrm{Lip}(X,E)$ satisfies $f \equiv 0$ on 
$\mathfrak{B} (T^{-1})$. Then $Tf \equiv 0$ on $\mathfrak{B} (T)$.
\end{lemma}

\begin{proof}
Suppose on the contrary that there exists $y_0 \in \mathfrak{B} (T)$ with $T f (y_0)  \neq 0$. By
Lemma~\ref{paco}, we can find $g \in \mathrm{Lip}(Y,F) $ with 
$\vc g(y_0) \vd = \vc g \vd_{\infty} > L(g)$
such that 
$$\vc g (y_0) \vd + \vc Tf (y_0) \vd = \vc (g + Tf) (y_0) \vd = \vc g + Tf \vd_{\infty} > L (g + Tf).$$

We see that 
\begin{eqnarray*}
\sup_{x \in \mathfrak{B} (T^{-1}) } \vc T^{-1} g (x) + f(x) \vd 
&=&\sup_{x \in \mathfrak{B} (T^{-1}) } \vc T^{-1} g (x) \vd \\ 
&\le&  \vc g \vd \\
&<& \vc (g + Tf) (y_0) \vd.
\end{eqnarray*}
On the other hand, if we put $ \mathrm{f} := (g + Tf) (y_0) $, since  $T^{-1}  \tilde{\mathrm{f}} \equiv 0$ 
on $\mathfrak{A} (T^{-1})$, there exists $n \in \mathbb{N}$ such that
\begin{eqnarray*}
\sup_{x \in \mathfrak{A} (T^{-1}) }  \vc T^{-1} g (x) + f(x) +  n T^{-1}  \tilde{\mathrm{f}} (x) \vd 
&<& \sup_{x \in \mathfrak{B} (T^{-1}) }  \vc T^{-1} g (x) + f(x) + n T^{-1}  \tilde{\mathrm{f}} (x) \vd \\
&<& (n+1) \vc \mathrm{f} \vd,
\end{eqnarray*}
so if we denote $k:=  T^{-1} \pl  g  + n \tilde{\mathrm{f}} \pr + f $, then we see that
$\vc k \vd_{\infty} < \vc Tk \vd$. Consequently,  $\vc k \vd_{\infty} < L (k)$ and there exists
$\mathrm{e} \in E$ with $T \tilde{\mathrm{e}} (y_0) \neq 0$  such that 
 $\vc k \pm  \tilde{\mathrm{e}} \vd_{\infty}
< L(k) = L(k \pm \tilde{\mathrm{e}})$. 

Also $L (Tk) = L (g + Tf) <\vc \mathrm{f} \vd$, so if we assume $n$ big enough, then $\vc Tk \vd = \vc Tk \vd_{\infty}$. 
 Therefore, $$
\vc  T \pl k \pm \tilde{\mathrm{e}} \pr \vd_{L} = \vc k \pm \tilde{\mathrm{e}} \vd_{L}=
L(k) = \vc T k \vd_{\infty} = \pl n+1 \pr \vc \mathrm{f} \vd .$$
This implies that 
$$ \vc \pl n+1 \pr  \mathrm{f} \pm  T \tilde{\mathrm{e}} (y_0)  \vd  = \vc  T \pl k \pm \tilde{\mathrm{e}} \pr (y_0)  \vd 
\le \pl n+1 \pr \vc \mathrm{f} \vd ,$$
which goes against Inequality~(\ref{strictly}).
\end{proof}

 Using Proposition~\ref{siegediesref},
we see that  the subspace $$\bxe := \tl f \in  \cxe : f \pl \as \pr \equiv 0 \tr$$ is isometrically isomorphic to $\bie$,
 via the restriction map. 
In the same way,  $$\axe := \tl f \in  \cxe : f \pl \bs \pr \equiv 0 \tr$$
and  $\aie$ are isometrically isomorphic.  Let denote by $I_{\bs} : \bie \ra \bxe$ and $I_{\as} : \aie \ra \axe$,
respectively, the corresponding natural isometries. 
In particular we can write in a natural way $$\cxe = \axe \oplus \bxe = \aie \oplus \bie ,$$ where this equality
has to be seen as a direct sum just  in the {\em linear} sense.

Next, let $R_{\bt} : \cyf \ra \bif$ be the operator sending each function to its restriction.

\begin{lemma}\label{parados}
The map  $$T_{\bh} := R_{\bt} \circ T \circ I_{\bs} : \bie \ra \bif $$ is a surjective linear isometry.
\end{lemma}

\begin{proof}
Notice first that if $f \in \bxe$ and $g \in \cxe$ satisfy $f \equiv g $ on $\bs$, then $\vc f \vd \le \vc g \vd$.

$T_{\bh}$ is linear and, by Lemma~\ref{1750ref}, it is easy to check that  it is surjective. We next see that it is 
an isometry.  Of course this is equivalent to show that $\vc R_{\bt} \circ T (f) \vd = \vc T (f) \vd$ 
for every $f \in \bxe$, and it is clear that $\vc R_{\bt} (T (f) ) \vd \le  \vc T (f) \vd$. 
Since $\vc R_{\bt} (T (f) ) \vd = \vc I_{\bt} \pl R_{\bt} (T (f) ) \pr \vd$, the fact that  
$\vc R_{\bt} (T (f) ) \vd <  \vc T (f) \vd$ is equivalent to that 
$$\vc T^{-1} \pl I_{\bt} \pl R_{\bt} \pl T \pl f \pr \pr \pr \pr \vd < \vc f \vd ,$$
which goes against the first comment in this proof.
\end{proof}

\begin{lemma}\label{altziturria}
$T_{\ah} := {I_{\at}}^{-1} \circ T \circ I_{\as} : \aie \ra \aif $ is standard.
\end{lemma}

\begin{proof}
Suppose that this is not the case. Since $\at = \ath \cup \bth$, we are in fact saying that $\ath \neq \emptyset$.
 
For $\mathrm{e} \in E$ with $\vc \mathrm{e} \vd =1$, we have 
$$ \vc T \pl   \tilde{\mathrm{e}} +   \chi_{\ash} \cdot \mathrm{e} \pr \vd =2 .$$
Notice  that both $T \tilde{\mathrm{e}}   \equiv 0$ and 
$T \pl \chi_{\as} \cdot \mathrm{e} \pr \equiv 0$ on $\ath$, 
so $T \pl \chi_{\bs} \cdot \mathrm{e} \pr \equiv 0$  on $\ath$. On the other hand, by Lemma~\ref{1750ref}, 
$c \pl T \pl \chi_{\ash} \cdot \mathrm{e} \pr \pr \subset \ath$, and consequently, since 
 $$ \vc T    \tilde{\mathrm{e}} \vd =  \vc T \pl    \chi_{\ash} \cdot \mathrm{e} \pr \vd =1 =  \vc T \pl   \tilde{\mathrm{e}} +   \chi_{\ash} \cdot \mathrm{e} \pr \vd_{\infty} ,$$  
 there are sequences  $(y_n)$ in $\ath$ and $(z_n) $ in $\bt$ with 
\begin{eqnarray*}
2 &=& \lim_{n \ra \infty} \frac{\vc  T \pl   \tilde{\mathrm{e}} +   \chi_{\ash} \cdot \mathrm{e} \pr   (y_n) -  T \pl   \tilde{\mathrm{e}} +   \chi_{\ash} \cdot \mathrm{e} \pr   (z_n)  \vd}
{d \pl y_n, z_n \pr} 
\\ &=& \lim_{n \ra \infty} \frac{\vc T \pl \chi_{\ash} \cdot \mathrm{e} \pr (y_n) -  T \pl \chi_{\bs} \cdot \mathrm{e} \pr  (z_n)  \vd}
{d \pl y_n, z_n \pr} \\
&=& \lim_{n \ra \infty} \frac{\vc  T \pl   \chi_{\bs} \cdot \mathrm{e} +   \chi_{\ash} \cdot \mathrm{e} \pr   (y_n) -  T \pl   \chi_{\bs} \cdot \mathrm{e} +   \chi_{\ash} \cdot \mathrm{e} \pr   (z_n)  \vd}
{d \pl y_n, z_n \pr} \\
&\le& \vc T \pl  \pl  \chi_{\bs} +   \chi_{\ash}  \pr \cdot \mathrm{e} \pr \vd \\
&\le& 1.
\end{eqnarray*}

We conclude that $\ath$ is empty.
\end{proof}

 It is easy to check that  $T_{\bh}$ satisfies Property {\bf P}, so by Theorem~\ref{marthahoyo}, it is standard. 
 We deduce the following result, which allows us to give the values on $\bt$ and on $\at$ of the images of all functions in $\cxe$ and $\axe$, respectively.
 
 \begin{corollary}\label{josecerecedo}
 There exist 
\begin{enumerate}
\item  $\hb \in \isob$   and  $\ha \in \isoa $, 
\item and maps $\jb : \bt \ra I(E,F)$ and $\ja: \at \ra I(E,F)$ 
constant on each $2$-component of $\bt$ and $\at$, respectively,
\end{enumerate}
 such that 
\begin{enumerate}
\item $
Tf(y)= \jb y (f(\hb (y)))
$
for all $f\in \cxe$ and $y\in \bt$, and
\item\label{facundo} $
Tf(y)= \ja y (f(\ha (y)))
$
for all $f\in \axe$ and $y\in \at$.
\end{enumerate}
\end{corollary}

\begin{lemma}\label{silvieta}
Let $y_0 \in \at$ and $A \subset \bs$ be such that $d(\ha (y_0), A) =1$. If $f \in \bxe$ satisfies
$f (A) \equiv \mathrm{e} \in E$, then $$T f (y_0) =  - \ja y_0 \pl \mathrm{e} \pr .$$
\end{lemma}
\begin{proof}
Notice first that, since $y_0 \in \at$, $ T \tilde{\mathrm{e}} (y_0) = 0$, and consequently, by Corollary~\ref{josecerecedo},  $  T \pl \chi_{\bs} \cdot \mathrm{e} \pr \pl y_0 \pr =  - T \pl \chi_{\as} \cdot \mathrm{e} \pr \pl y_0 \pr = - \ja y_0 \pl \mathrm{e} \pr  $.

Next we  prove the result through several steps. We denote $\mathrm{a} :=  \ja y_0 \pl \mathrm{e} \pr$ for short.

\begin{step}\label{zelada}
Assume that $\vc \mathrm{e} \vd = 1 =\vc f \vd$. 
\end{step}

Consider $k' \in \cx$ defined by $k ' (x) := \max \tl 0, 1 - d(x, \ha (y_0))  \tr $
for every $x \in X$, and $k \in \axe$ defined by $k := - k' \cdot \mathrm{e}$. It is easy to see that $(k + f ) (\ha (y_0)) = - \mathrm{e}$ and that $(k + f ) (x) =  \mathrm{e}$ for every $x \in A$. As a consequence, 
$\vc k + f \vd =2$. 

Suppose now that $T f (y_0) = \mathrm{f} \neq -   \mathrm{a}  $. By  Corollary~\ref{josecerecedo}, $Tk (y_0) = - \mathrm{a}$
and, since $\vc  f  \vd_{\infty}  =1$, we can take $M <2$ such that 
$$
\vc T (k + f) (y_0) \vd =  \vc  \mathrm{-a} + \mathrm{f}  \vd 
< M.
$$
Consequently there exists $0<r<1$ such that $\vc T \pl k + f     \pr (y) \vd < M$ for every 
$y \in B(y_0, r)$. On the other hand, for $y \in \at$ with $d(y, y_0 ) \ge r$, 
\begin{eqnarray*}
\vc T k (y) \vd &=& \vc T (- k' \cdot \mathrm{e} )   (y) \vd \\
&=& \vc    \max \tl 0, 1 - d(\ha (y), \ha (y_0))  \tr  \cdot  \ja y (\mathrm{e})   \vd \\
&\le& 1- r,
\end{eqnarray*}
so $\vc T (k + f) (y) \vd \le 2-r$. Since $\vc T (k + f) (y) \vd = \vc T f (y) \vd \le 1$ for every $y \in \bt$,
 we deduce that $\vc T (k + f)  \vd_{\infty} < 2 = \vc T (k + f)  \vd$. Let $M>0$ with $M + \vc T (k + f)  \vd_{\infty} <2$, and  
$y \in \bt$ such that $\hb (y) \in A$ and $d(\hb (y), \ha (y_0)) < 1+M/2$. Define $\mathrm{b} := M T \tilde{\mathrm{e}} (y) \in F$. 
 By Corollary~\ref{josecerecedo},   $T^{-1} \tilde{\mathrm{b}} (\hb (y)) = M \mathrm{e}$, and consequently 
 \begin{eqnarray*}
2 d (\hb (y), \ha (y_0))    &<& 2 + M  \\
&=& \vc \pl k + f + T^{-1} \tilde{\mathrm{b}} \pr \pl \hb \pl y \pr \pr - \pl k + f + T^{-1} \tilde{\mathrm{b}} \pr \pl \ha \pl y_0 \pr \pr \vd ,
\end{eqnarray*}  
 against the fact that $\vc   k + f + T^{-1} \tilde{\mathrm{b}} \vd =2$.

\begin{step}\label{sinapsis}
Assume that $\vc \mathrm{e} \vd = 1 =\vc f \vd_{\infty}$.  
\end{step}

It is easy to check that  if $n \ge L (  f)$, 
then $\vc    n \chi_{\bs} \cdot \mathrm{e}  \vd_{\infty} = \vc   n \chi_{\bs} \cdot \mathrm{e}  \vd = n$  and that 
$$\vc  f +  n \chi_{\bs} \cdot \mathrm{e}  \vd_{\infty} = \vc  f +  n \chi_{\bs} \cdot \mathrm{e}  \vd = n+1 .$$ 
Using Step~\ref{zelada}, $T (   n \chi_{\bs} \cdot \mathrm{e} ) (y_0) = - n \mathrm{a}$ and
$T ( f +  n \chi_{\bs} \cdot \mathrm{e} ) (y_0) = - (n +1)  \mathrm{a}$. The conclusion is easy.

\begin{step}\label{aurillaque}
Assume that $\mathrm{e} =0$.
\end{step}

Of course we must prove that $Tf (y_0) =0$. Fix $\mathrm{d} \in E$ with norm $1$. Consider $m \in \bxe$ defined by 
$
m (x):= \max \{0,1-d(x,A)\} \cdot \mathrm{d}
$
for each $x\in X$. 
We easily check that $\vc   m \vd_{\infty}= 1 = \vc \mathrm{d} \vd$,
and  if we assume that $\vc   f \vd \le 1$, then $\vc f(x) \vd \le d(x, A)$ for every $x$. 
As in  the proof of 
Lemma~\ref{paco}, we see that 
$\vc   m +  f  \vd_{\infty}= 1 = \vc \mathrm{d} \vd$.
The conclusion follows immediately from Step~\ref{sinapsis}. 

\medskip

The rest of the proof  is easy.
\end{proof}

\begin{corollary}\label{indoxina}
Suppose that $A_1, A_2 \subset \bs$ satisfy $d(\ha (y_0), A_i ) =1$ for $i=1,2$. Then $d(A_1, A_2) = 0$.
\end{corollary}

\begin{proof}
Just assume that $d(A_1, A_2) > 0$ and apply  Lemma~\ref{silvieta} 
to any  $  f \in \bxe$  
such that $  f   (A_i) \equiv (-1)^i \mathrm{e} \neq 0$ for $i=1, 2$. This leads to two different values for  
$Tf (y_0)$.
\end{proof}

\begin{corollary}\label{quena}
Let $y_0 \in  \at$. Then there exists exactly one point $\varphi (y_0) $ in $ \bt$
such that $d \pl \hb \pl \varphi (y_0) \pr, \ha \pl y_0 \pr  \pr = 1$. Also,  
$$T f (y_0) =  -   \ja y_0 (f \pl \hb \pl \varphi \pl  y_0 \pr \pr  \pr ) $$
for every $f \in \bxe$.
\end{corollary}

\begin{proof}
By Lemma~\ref{gatosinchicle} and Corollary~\ref{indoxina}, we deduce that if $(x_n)$ is a sequence in $X$ such that $d (\ha (y_0), x_n ) \le 1 + 1/n$ 
for each $n \in \mathbb{N}$, then it is a Cauchy sequence, so there is a limit $x_0$ in $ \widehat{X}$, which necessarily belongs to $\widehat{\bs}$. Obviously the point $x_0$ does not depend on the sequence we take.

We next assume that $\widehat{X}$ is not complete
 and prove that
$x_0 \in \bs$. 
 If this is not the case,
for each $n \in \mathbb{N}$, let 
$$
f_n (x) := \max \tl 0,1 -   d \pl x,  B \pl x_0 , 1/n \pr \pr   \tr 
$$
for all $x \in   X$. It is clear that each $f_n$ belongs to  
$\cx$. Since $\mathrm{Lip}(X,\widehat{E})$ is complete, if we take $(\mathrm{e}_n)$ in $E$ as in Lemma~\ref{ester}, then 
$f:= \sum_{n=1}^{\infty} f_n \cdot \widetilde{\mathrm{e}_n}$ belongs to 
$\mathrm{Lip}(X,\widehat{E})$, and since all values are taken in $E$, to  $\mathrm{Lip}(X,E)$, and indeed to $\bxe$.
Thus, since $f = \lim_{k \ra \infty } \sum_{n=1}^k f_n \cdot \widetilde{\mathrm{e}_n}$, we deduce 
from Lemma~\ref{silvieta} 
that  
\begin{eqnarray*}
 Tf (y_0) &=& \lim_{k \ra \infty}  \sum_{n=1}^k T \pl f_n \cdot  \widetilde{\mathrm{e}_n} (y_0) \pr 
\\  &=&  -  \lim_{k \ra \infty} \sum_{n=1}^{k}   \ja y_0  \pl \mathrm{e}_n  \pr \\
 &=& -     \sum_{n=1}^{\infty} \ja y_0 \pl \mathrm{e}_n  \pr ,
\end{eqnarray*}
  which belongs to $\widehat{F} \setminus  F$. This is absurd.
  
  If we define $\varphi (y_0) := {h_{\bh}}^{-1} (x_0) \in \bt$, then we are done.
\end{proof} 

\begin{proposition}\label{nohayna}
For every $y \in \at$, $\ja y = - \jb \varphi (y)$.
\end{proposition}

\begin{proof}
Fix $y \in \at$ and $\mathrm{e} \in E$, and let $\mathrm{f} := T \tilde{\mathrm{e}} (\varphi (y))$. Then 
$T^{-1} \tilde{\mathrm{f}} (\hb (\varphi (y))) = \mathrm{e}$ and $T^{-1} \tilde{\mathrm{f}} \equiv 0$ on $\as$.
We conclude from Corollary~\ref{josecerecedo} that $\jb \varphi (y) (\mathrm{e}) = \mathrm{f}$,
 and from Corollary~\ref{quena} that $- \ja y (\mathrm{e}) = \mathrm{f}$.
\end{proof}

Next result follows now easily from Corollaries~\ref{josecerecedo} and~\ref{quena}, and Proposition~\ref{nohayna}.

\begin{corollary}\label{errabia}
For $y \in \at$ and $f \in \cxe$, 
\begin{eqnarray*}
Tf (y) &=& -\ja y (f (\hb (\varphi (y)))) + \ja y (f (\ha (y))) \\
&=& \jb \varphi (y)  (f (\hb (\varphi (y)))) - \jb \varphi (y) (f (\ha (y))) .
\end{eqnarray*}
\end{corollary}

\begin{corollary}\label{kamisabieja}
Let $y_0 \in \at$. If $y \in \bt$ is such that   $  d(y, \varphi (y_0 )) \ge 2$, then
$$d (y, y_0) \ge 2 .$$
\end{corollary}

\begin{proof}
Let $\mathrm{e}_1, \mathrm{e}_2 \in E$ be vectors with norm $1$   and such that $\jb \varphi(y_0) (\mathrm{e}_1) = \jb y (\mathrm{e}_2)  $. Define $f:= f_1 - f_2 \in \bxe$, where
 $$f_1 (x) := \max \tl 0, 1- d(x, \hb \pl \varphi (y_0) \pr ) \tr \cdot \mathrm{e}_1$$
and $$ f_2 (x) := \max \tl 0, 1- d(x, \hb \pl y \pr ) \tr \cdot \mathrm{e}_2 $$
 for every $x \in X$. Obviously, $\vc f_1 \vd = 1 = \vc f_2 \vd$, so to show that $\vc f \vd =1$,
 it is enough to see that, if $d(x, \hb \pl \varphi (y_0) \pr ) , d(z, \hb \pl y \pr ) <1$, then 
 $\vc f_1 (x) \vd + \vc f_2 (z) \vd \le d(x, z)$. Taking into account that
 $$2 \le  d(\hb \pl y \pr , \hb \pl \varphi (y_0 ) \pr ) \le d(z, \hb \pl y \pr ) + d(x,z ) + d(x, \hb \pl \varphi (y_0) \pr ) , $$
 it follows that 
 \begin{eqnarray*}
 \vc f_1 (x) \vd + \vc f_2 (z) \vd &=& \pl 1- d(x, \hb \pl \varphi (y_0) \pr ) \pr + \pl 1 - d(z, \hb \pl y \pr ) \pr \\
 &\le& d(x, z).
 \end{eqnarray*}
 
On the other hand,  by Corollary~\ref{errabia}, $ Tf (y_0) =  \jb \varphi (y_0) (\mathrm{e}_1) = Tf (\varphi (y_0))$, 
and by the way we have taken $\mathrm{e}_1$ and $ \mathrm{e}_2 $, we have
 $ Tf (y) = - \jb y  (\mathrm{e}_2) = - \jb \varphi(y_0) (\mathrm{e}_1) $.
 We conclude that, since $\vc Tf \vd = 1$, $$2 = \vc Tf (y) - Tf (y_0) \vd \le d(y, y_0). $$ 
\end{proof}

\begin{corollary}\label{nitanmal}
Let $y_0 \in \at$. Given $y \in \bt$, if  $0 \le d(y, \varphi (y_0 )) <1$, then
$$d (y, y_0) = 1 + d (y, \varphi (y_0 )) ,$$
and if $1 \le d(y, \varphi (y_0 )) <2$, then
$$d (y, y_0) \ge 2 .$$
\end{corollary}

\begin{proof}
Fix $\mathrm{e} \in E$ with norm $1$ and let 
$f (x) :=  \min \tl 1,  d(x , \hb (\varphi (y_0)) ) \tr \cdot \mathrm{e} $ for every $x \in X$.
Let $y \in \bt$ with $ d(y, \varphi (y_0 )) <2$. Taking into account that 
 $\vc f \vd =1$,  Corollaries~\ref{josecerecedo} and~\ref{errabia} give 
\begin{eqnarray*}
d(y, y_0)  &\ge& \vc Tf (y) - T f (y_0) \vd\\
 &=& \vc \min \tl 1,  d(\hb (y) , \hb (\varphi (y_0)) ) \tr \cdot \jb y (\mathrm{e}) + \jb y (\mathrm{e}) \vd \\
  &\ge& \min \tl 2, d( \hb (y) , \hb (\varphi (y_0)) ) + 1 \tr.
\end{eqnarray*}
The conclusion is immediate.
\end{proof}

\begin{corollary}\label{tillohablador}
If $y_1, y_2 \in \at$ satisfy $\varphi (y_1) \neq \varphi (y_2)$, then 
$d (y_1, y_2 ) \ge 2  $.
\end{corollary}

\begin{proof}
Suppose that $M:= d(y_1, y_2) /2 <1$, so by Corollary~\ref{josecerecedo} $\ja y_1 = \ja y_2$.
Put $ N:=  d(\hb (\varphi (y_1)) , \hb (\varphi (y_2))) $ and,  
 for a fixed $\mathrm{e} \in E$ with norm $1$, let 
$$f (x ) := \max \tl - M, M - d( x , \ha (y_1) ) \tr \cdot \mathrm{e} $$
and  
$$g (x) := \max \tl 0 , N - d(x , \hb (\varphi (y_1)) ) \tr \cdot \mathrm{e} $$
for every $x \in X$.
If we take  
 $A >0$ such that $M + AN <1$, then 
$k := \chi_{\as} f  - A \chi_{\bs} g $ has norm $1$. Also, by Corollary~\ref{errabia}, $Tk (y_1) = M \ja y_1 (\mathrm{e}) +  A N \ja y_1 (\mathrm{e})$ and 
$Tk (y_2) = - M \ja y_2 (\mathrm{e}) = - M \ja y_1 (\mathrm{e})$. Consequently 
$$ \vc \frac{Tk (y_1) - Tk(y_2)}{d(y_1, y_2)} \vd >1,$$ which is impossible. 
\end{proof}

\begin{proof}[Proof of Theorem~\ref{nonelk}]
 Corollaries~\ref{kamisabieja}, ~\ref{nitanmal} and \ref{tillohablador} show that $Y$ is of type $\mathbf{A}$. 
We consider the associated
purely nonstandard map $S_{\varphi} : \cyf \ra \cyf$ and see that, given $\mathrm{e} \in E$, $\mathrm{e} \neq 0$,
 the composition $S_{\varphi} \circ T : \cxe \ra \cyf$
satisfies $S_{\varphi} \circ T \pl \tilde{\mathrm{e}} \pr = \jb \varphi (y) (\mathrm{e}) \neq 0$ if $y \in \at$ and
$S_{\varphi} \circ T \pl \tilde{\mathrm{e}} \pr = \jb   y (\mathrm{e}) \neq 0$ if $y \in \bt$, that is, the composition satisfies Property {\bf P}.

This implies that $S_{\varphi} \circ T$ is standard.  Since $S_{\varphi} = {S_{\varphi}}^{-1}$, we have that
$T = S_{\varphi} \circ \pl S_{\varphi} \circ T \pr$, and we are done.
\end{proof}

\begin{remark}
It is easy to check that, if  $h: Y \ra X$
and $J : Y \ra I(E, F)$ are the 
associated maps to $S_{\varphi} \circ T$ in the proof of Theorem~\ref{nonelk}, then
 $h \equiv \ha$ on $\at$ and $h \equiv \hb$ on $\bt$. In the same way,
$J \equiv - \ja$ on $\at$ and $J \equiv  \jb$ on $\bt$. Finally, it is also apparent that, 
given $A \in \mathtt{Comp}_2 (\bs)$, $A = A' \cap \bs$, where $A' \in \mathtt{Comp}_2 (X)$, and that
a similar fact does not necessarily hold for the elements in 
 $ \mathtt{Comp}_2 (\as)$.
\end{remark}

\begin{proof}[Proof of Proposition~\ref{pesafrank}]
Suppose that $\pl X, d^{\alpha}\pr$ is of type $\mathbf{A}$. Then  
$$d^{\alpha} (z, y) = 1 + d^{\alpha} (z, \varphi (y)) $$ whenever 
$y \in \at$ and $z \in \bt$ satisfy  $0 < d^{\alpha} (z, \varphi (y)) <1$. 
In such case,  $d^{\alpha} (z, y) < d  (z, y) $ and $d (z, \varphi (y )) < d^{\alpha} (z, \varphi (y))  $, 
and this implies 
$$d (z, y) > 1 + d (z, \varphi (y)) = d(y , \varphi (y) )  + d (z, \varphi (y)) , $$
which is impossible. We deduce that, if $d  (z, \varphi (y )) <1$, then $z = \varphi (y)$. 
The rest of the proof is easy.
\end{proof}

\section{Acknowledgements}
The authors wish  to thank the referee for his/her  valuable remarks, which helped them to improve
the first version of the paper.

\end{document}